\documentclass[12pt]{article}
\usepackage{amsfonts}

\newtheorem{prop}{Proposition}[section]
\newtheorem{cor}[prop]{Corollary}
\newtheorem{unnumber}{}

\newtheorem{prepf}[unnumber]{Proof}
\newenvironment{pf}{\prepf\rm}{\endprepf}
\newcommand{\qed}{\hfill$\Box$}

\begin{document}

\title{Some measures of finite groups related to permutation bases}
\author{Peter J. Cameron}
\date{School of Mathematical Sciences\\Queen Mary, University of London\\
Mile End Road\\London E1 4NS\\U.K.\\\texttt{p.j.cameron@qmul.ac.uk}}
\maketitle

\begin{quote}
\textbf{Note:} This paper is more than ten years old and has been rejected
by a couple of journals. I am posting it on the arXiv (after minimal editing)
since I think there is still some interest in these questions.
\end{quote}

\begin{abstract}
I define three ``measures'' of the complicatedness of a finite group in terms
of bases in permutation representations of the group, and consider their
relationships to other measures.
\end{abstract}

\section{Introduction and history}

The purpose of this paper is to define and investigate three functions which
in some sense ``measure'' the complicatedness of a finite group. The functions
are defined in minimax fashion in terms of base size in permutation
representations of the group. The values of these measures for symmetric
groups will be given.

I begin with a brief account of some similar measures, which are related to
the new measures in various ways.

If $f$ is a function from (isomorphism classes of) finite groups to natural
numbers, let $f'$ be the function defined by
\[f'(G)=\max\{f(H):H\le G\}.\]
The function $f'$ is necessarily monotonic: that is, if $H\le G$, then
$f'(H)\le f'(G)$. Conversely, if $f$ is monotonic, then $f=f'$; so every
monotonic function is of this form.

An obvious function is $d(G)$, the minimum number of generators of $G$.
Now is is well-known that $d(S_n)=2$; so this invariant is not very sensitive.

However, $d'$ is more interesting. Babai~\cite{babai} raised the question of
determining $d'(S_n)$. The motivation was computational group theory: if we
are given an arbitrary (possibly very large) subset of $S_n$, and want to
know the size of the subgroup it generates, then we know that we can replace
the given set by a set of size at most $d'(S_n)$ generating the same subgroup.
This suggests that we also want, if possible, an algorithmic method for
finding such a minimum-size generating set.

Babai considered the related function $l$, where $l(G)$ is the length of the
longest chain of subgroups of $G$. Since, obviously, $d(G)\le l(G)$, and
$l$ is monotonic, (that is, $l(G)=l'(G)$), we see that $d'(G)\le l(G)$
for any group $G$. Babai proved that
\[d'(S_n)\le l(S_n)\le 2n-1.\]

The exact value of $S_n$ was established by Cameron, Solomon and
Turull~\cite{cst}:
\[l(S_n)=\left\lceil\frac{3n}{2}\right\rceil-b(n)-1,\]
where $b(n)$ is the number of $1$s in the base~$2$ representation of $n$.

Jerrum~\cite{jerrum} showed that $d'(S_n)\le n-1$. His proof was
algorithmic: given an arbitrary set of permutations, a set of at most $n-1$
permutations generating the same subgroup can be found efficiently.

The exact value of $d'(S_n)$ was established by McIver and Neumann~\cite{mn}:
\[d'(S_n)=\left\lfloor\frac{n}{2}\right\rfloor\hbox{ for }n>3.\]

Whiston~\cite{whiston} considered the invariant $\mu(G)$, the maximal size
of an independent generating set for $G$, where a set is \emph{independent}
if none of its elements lies in the subgroup generated by the others. Since
any independent set is an independent generating set for the subgroup it
generates, $\mu'(G)$ is the maximum size of an independent subset of $G$.
The parameter $\mu(G)$ appears in the work of Diaconis and
Saloff-Coste~\cite{ds} in the rate of convergence of the product replacement
algorithm for finding a random element of a finite group.

Whiston showed that
\[\mu(S_n)=\mu'(S_n)=n-1.\]
However, he observed that there are groups with $\mu'(G)>\mu(G)$.

Cameron and Cara~\cite{cc} found all independent generating sets of size
$n-1$ in $S_n$.

\section{The base measures}

Let $G$ be a permutation group on $\Omega$. A \emph{base} for $G$ is a
sequence of points of $\Omega$ whose pointwise stabiliser is the identity.
(Treating a base as a sequence rather than a set fits in with the use of
bases in computational group theory, where the elements of a base are chosen
in order.) A base is called \emph{irredundant} if no point is fixed by the
pointwise stabiliser of its predecessors; it is \emph{minimal} if no point
is fixed by the pointwise stabiliser of the other points in the sequence.

It is computationally a simple matter to choose an irredundant base; simply
choose each base point to be a point moved by the stabiliser of the points
previously chosen. It is less straightforward to choose a minimal base.
Note that a base is minimal if and only if every re-ordering of it is
irredundant.

Now, given a finite group $G$, we define three numbers $b_1(G)$, $b_2(G)$,
$b_3(G)$ as follows. In each case, the maximum is taken over all permutation
representations of $G$ (not necessarily faithful).
\begin{itemize}
\item $b_1(G)$ is the maximum, over all representations, of the maximum size
of an irredundant base;
\item $b_2(G)$ is the maximum, over all representations, of the maximum size
of a minimal base;
\item $b_3(G)$ is the maximum, over all representations, of the minimum
base size.
\end{itemize}
Clearly we have:

\begin{prop}
$b_3(G)\le b_2(G)\le b_1(G)$.
\end{prop}

These inequalities can be strict. The group $G=\mathrm{PSL}(2,7)$ has
$b_1(G)=5$, $b_2(G)=4$, and $b_3(G)=3$.

Now $b_1(G)$ is a parameter we have seen before!

\begin{prop}
$b_1(G)=l(G)$.
\end{prop}

\begin{pf}
Given an irredundant base of size $b_1(G)$, the stabilisers of its initial
subsequences form a properly descending chain of subgroups of length 
$b_1(G)$. So $b_1(G)\le l(G)$.

Conversely, let
\[G=G_0>G_1>\cdots>G_l=1\]
be a chain of subgroups of length $l=l(G)$. Consider the action on the union
of the coset spaces of the subgroups $G_i$, and let $\alpha_i$ be the point
$G_i$ of the coset space $(G:G_i)$. Then $(\alpha_1,\ldots,\alpha_l)$ is
an irredundant base. So $l(G)\le b_1(G)$.\qed
\end{pf}

We will see in the next section a connection between $b_2(G)$ and $\mu(G)$.
I know much less about $b_3(G)$. One observation is the following:

\begin{prop}
Let $G$ be a non-abelian finite simple group. Then $b_3(G)$ can be calculated
by considering only the primitive permutation representations of $G$.
\end{prop}

\begin{pf}
Given any permutation representation of $G$, we can discard fixed points, so
that $G$ acts faithfully on each orbit. Now let $b_3^*(G)$ be the maximum
of the minimum base sizes over all transitive representations of $G$, and
suppose that there is an intransitive representation with minimum base size
greater than $b_3^*(G)$. Now there exist at most $b_3^*(G)$ points in an orbit
whose stabiliser acts trivially on that orbit, and hence is trivial (since
the action on the orbit is faithful), contrary to assumption.

Now let $b_3^+(G)$ be the maximum of the minimum base sizes over all primitive
representations of $G$, and suppose that there is a transitive but
imprimitive representation with base size greater than $b_3^+(G)$. There are
at most $b_3^+(G)$ maximal blocks whose stabiliser acts trivially on the
block system, and hence is trivial (since again the action is faithful),
contrary to assumption.\qed
\end{pf}

This proposition does not hold for $b_2(G)$. For the group
$G=\mathrm{PSL}(2,7)\cong\mathrm{PSL}(3,2)$, a minimal base in any transitive
representation has size at most $3$. However, in the action on the points
and lines of the projective plane of order~$2$, there is a minimal base of
size~$4$, consisting of two points and two lines such that each point lies
on one of the lines and each line passes through one of the points.

\section{Boolean semilattices}

The subgroups of the group $G$ form a lattice $L(G)$, with the operations
$H\wedge K=H\cap K$ and $H\vee K=\langle H,K\rangle$. A 
\emph{meet-semilattice} of $L(G)$ is a collection of subgroups containing $G$
and closed under $\wedge$, while a \emph{join-semilattice} is a collection
of subgroups containing the trivial group $1$ and closed under $\vee$.

The \emph{Boolean lattice} $B(n)$ is the lattice of subsets of an $n$-set.

\begin{prop}
Let $G$ be a finite group. Then $B(n)$ is embeddable as a meet-semilattice
in $L(G)$ if and only if it is embeddable as a join-semilattice.
\end{prop}

\begin{pf}
Suppose first that $B(n)$ is a join-semilattice of $L(G)$. Let
$N=\{1,\ldots,n\}$. Then, for every subset $I$ of $N$, there is a subgroup
$H_I$ of $G$, and $H_{I\cup J}=\langle H_i,H_j\rangle$ for any two subsets
$I$ and $J$. Moreover, all these subgroups are distinct. In particular,
$H_i\not\le H_{N\setminus\{i\}}$ for all $i$ (where $H_i$ is shorthand for
$H_{\{i\}}$); else
\[H_N=\langle H_i,H_{N\setminus\{i\}}\rangle=H_{N\setminus\{i\}},\]
contrary to assumption.

Let $K_i=H_{N\setminus\{i\}}$, and for any $I\subseteq N$, put
\[K_I=\bigcap_{i\in I}K_i,\]
with the convention that $K_\emptyset=G$. We claim that all the subgroups
$K_i$ are distinct. Suppose that two of them are equal, say $K_I=K_J$.
By interchanging $I$ and $J$ if necessary, we may assume that there exists
$i\in I\setminus J$. But then $H_i\le K_J$ while $H_i\not\le K_I$, a
contradiction.

Now it is clear that $K_J\cap K_J=K_{I\cap J}$, so we have an embedding of
$B(n)$ as a meet-semilattice (where we have reversed the order-isomoprhism
to simplify the notation).

The reverse implication is proved by an almost identical argument.\qed
\end{pf}

Note that the conditions of the proposition are not equivalent to
embeddability of $B(n)$ as a lattice. For example, if $G$ is the quaternion
group of order $8$, then $B(2)$ is embeddable as both a meet-semilattice
and a join-semilattice but not as a lattice.

\begin{prop}
Let $G$ be a finite group.
\begin{enumerate}
\item The largest $n$ such that $B(n)$ is embeddable as a join-semilattice
of $L(G)$ is $\mu'(G)$.
\item The largest $n$ such that $B(n)$ is embeddable as a meet-semilattice
of $L(G)$ in such a way that the minimal element is a normal subgroup of $G$
is $b_2(G)$.
\end{enumerate}
\end{prop}

\begin{pf}
(a) Let $\{g_1,\ldots,g_n\}$ be an independent set in $G$, where $n=\mu'(G)$.
Let $N=\{1,\ldots,n\}$. Then the subgroups $H_I=\langle g_i:i\in I\rangle$
form a join-semilattice of $G$ isomorphic to $B(n)$.

Conversely, suppose that we have a join-semilattice given by the subgroups
$H_I$ for $I\subseteq N$. Choose $g_i\notin H_{N\setminus\{i\}}$. Then
clearly the elements $g_1,\ldots,g_n$ are independent.

\smallskip

(b) Let $(\alpha_1,\ldots,\alpha_n)$ be a minimal base for $G$ in some
permutation representation. Let $K_i$ be the stabiliser of $\alpha_i$, and
$K_I=\bigcap_{i\in I}K_i$ for $I\subseteq N=\{1,\ldots,n\}$. Then the 
subgroups $K_I$ form a meet-semilattice of $G$. The subgroup$ K_N$ is the
kernel of the permutation representation (by definition of a base), and so
is a normal subgroup of $G$.

Conversely, suppose that we have a meet-semilattice given by the subgroups
$K_I$ for $I\subseteq N$, such that $K_N$ is a normal subgroup fo $G$;
the notation is chosen so that $K_I\cap K_J=K_{I\cup J}$. Now consider the
permutation representation on the union of the coset spaces of the subgroups
$K_i=K_{\{i\}}$. Since the intersection of all these subgroups is normal,
it is the kernel of the representation, and the points corresponding to
the subgroups $K_i$ form a base. It is minimal, since the intersection
of fewer than $n$ of the subgroups is not equal to $K_N$.\qed
\end{pf}

\begin{cor}
$b_2(G)\le\mu'(G)$ for any group $G$.
\end{cor}

I do not know a group where the inequality is strict. Resolving this is
equivalent to the following question. Let $n$ be maximal such that $B(n)$
is embeddable as a meet-semilattice of $L(G)$. Is there an embedding of $B(n)$
for which the minimal element is normal?

From the above results, we conclude:

\begin{cor}
$b_2(S_n)=b_3(S_n)=n-1$.
\end{cor}

\begin{pf}
We have
\[n-1\le b_3(S_n)\le b_2(S_n)\le \mu'(S_n)=n-1,\]
wwhere the first inequality holds because any base in the natural
representation has size $n-1$; the second is trivial; the third comes from
the preceding Corollary; and the equality is Whiston's theorem.\qed
\end{pf}

\end{document}